\documentclass[a4paper,11pt,twoside,english]{article}
\usepackage{babel}
\usepackage{latexsym,amsfonts}
\usepackage{amsthm}
\usepackage{amsmath}
\usepackage{mathtools}
\usepackage{paralist}
\usepackage{todonotes}
\usepackage{tikz}

\usepackage{lastpage}
\usepackage{fancyhdr}
\fancypagestyle{plain}{

  \fancyhead[L,C]{}
  \fancyhead[R]{}
  \fancyfoot[R,L]{}
  \fancyfoot[C]{\thepage \ of \pageref{LastPage}}
  }

\title{A remark on nearness spaces}
\author{Jan-David Hardtke}
\date{}

\DeclareMathOperator{\interior}{int}

\providecommand{\sm}{\setminus}
\providecommand{\ssq}{\subseteq}
\providecommand{\id}{\ensuremath{\mathrm{id}}}
\providecommand{\A}{\ensuremath{\mathcal{A}}}
\providecommand{\B}{\ensuremath{\mathcal{B}}}

\providecommand{\circled}[1]{\tikz[baseline=(char.base)]{\node[shape=circle,draw,inner sep=1.5pt] (char) {#1};}}

\providecommand{\keywords}[1]{
  {\let\thefootnote=\relax
  \footnote{{\em Keywords}: #1}}
  \addtocounter{footnote}{-1}
  }

\providecommand{\AMS}[1]{
  {\let\thefootnote=\relax
  \footnote{{\em AMS Subject Classification} (2010): #1}}
  \addtocounter{footnote}{-1}
  }

\providecommand{\address}{
  {\sc \noindent Department of Mathematics \\
  Universit\"at Leipzig\\
  Augustusplatz 10, 04109 Leipzig \\
  Germany \\}
}

\DeclarePairedDelimiter{\set}{\lbrace}{\rbrace}

\theoremstyle{definition}

\theoremstyle{remark}

\newtheorem*{claim*}{Claim}

\newtheorem*{fact*}{Fact}
\theoremstyle{plain}

\newtheorem*{lemma*}{Lemma}

\newtheorem*{theorem*}{Theorem}

\newenvironment{Proof}[1][\proofname]{\begin{proof}[#1] \setlength{\parindent}{0pt}}{\end{proof}}
\newenvironment{Abstract}{\centering\begin{minipage}{0.8\textwidth} \noindent \small {\sc Abstract.}}{\end{minipage}\par}

\usepackage{color}
\definecolor{darkgreen}{rgb}{0,0.5,0}

\numberwithin{equation}{section}
\newtagform{colored}[\color{blue}]{\color{blue}(}{\color{blue})}
\usetagform{colored}

\usepackage[colorlinks,linkcolor=blue,citecolor=red,urlcolor=darkgreen]{hyperref}
\providecommand{\email}{{\it E-mail address:} \href{mailto:hardtke@math.uni-leipzig.de}{\tt hardtke@math.uni-leipzig.de}}

\usepackage{amsrefs}

\begin{document}

\maketitle

\begin{Abstract}
We give a proof of the well-known fact that the category of nearness spaces is bireflective in the
category of merotopic spaces which uses Zorn's Lemma instead of the usual construction by transfinite
induction.
\end{Abstract}
\AMS{54E17; 54B30}

\ \par
First let us introduce some notation and definitions. If $X$ is any set we denote by $\mathcal{P}(X)$ the power-set of $X$ 
and by $c(X)$ the set of all non-empty coverings of $X$. For $\A, \B\in c(X)$ we say that $\A$ is a refinement of $\B$ 
(or $\A$ refines $\B$), denoted by $\A\prec\B$, if for every $A\in \A$ there is $B\in \B$ such that $A\ssq B$. We further 
put $\A\wedge\B=\set*{A\cap B:A\in \A, B\in \B}$, which is obviously again a covering of $X$. For any unexplained notions from category theory which we will use in the sequel the reader is referred to \cite{preuss1} and \cite{preuss2}.\par
Now recall that a merotopic space is a pair $(X,\mu)$, where $\mu\ssq c(X)$ is non-empty and such that the following holds:
\begin{enumerate}[(i)]
\item $\A\in \mu, \B\in c(X) \ \mathrm{and} \ \A\prec\B \ \Rightarrow \ \B\in \mu$
\item $\A, \B\in \mu \ \Rightarrow \A\wedge\B\in \mu$
\end{enumerate}
Then $\mu$ is sometimes called a merotopic structure on $X$ and the elements of $\mu$ are called uniform coverings.\par
A map $f:X \rightarrow Y$ between merotopic spaces $(X,\mu)$ and $(Y,\nu)$ is called uniformly continuous 
(with respect to $\mu$ and $\nu$) if $f^{-1}[\A]\in \mu$ for every $\A\in \nu$, where $f^{-1}[\A]:=\set*{f^{-1}(A):A\in \A}$.\par
Clearly, the merotopic spaces together with the uniformly continuous maps as morphisms form a concrete category over the category 
of sets, i.\,e. a construct, which will be denoted by $\mathbf{Mer}$.\par
For a merotopic space $(X,\mu)$ and a subset $A\ssq X$ we define the interior of $A$ with respect to $\mu$ as $\interior_{\mu}(A)=
\set*{x\in X:\set*{A,X\sm\set{x}}\in \mu}$.\par
It is well-known and easily checked that the following assertions hold 
(cf. \cite{preuss1}*{Proposition 3.2.2.3}, where merotopic spaces are called semi-nearness spaces):
\begin{enumerate}[(a)]
\item $\interior_{\mu}(A)\ssq A \ \ \forall A\ssq X$
\item $\interior_{\mu}(X)=X, \ \interior_{\mu}(\emptyset)=\emptyset$
\item $A\ssq B\ssq X \ \Rightarrow \ \interior_{\mu}(A)\ssq \interior_{\mu}(B)$
\item $\interior_{\mu}(A\cap B)=\interior_{\mu}(A)\cap\interior_{\mu}(B) \ \ \forall A, B\ssq X$
\end{enumerate}\par
A merotopic space $(X,\mu)$ is called a nearness space if $\set*{\interior_{\mu}(A):A\in \A}$ belongs to $\mu$ for every $\A\in \mu$.
Then $\mu$ is also called a nearness structure on $X$. The nearness spaces (and uniformly continuous maps) induce a full subconstruct of 
$\mathbf{Mer}$ which will be denoted by $\mathbf{Near}$.\par
The merotopic spaces were originally introduced by Kat\v{e}tov in \cite{katetov}, though not in the formulation above
but in an equivalent version using so called micromeric collections of subsets of $X$ instead of uniform coverings. One can also use a concept of ``near'' collections of subsets of $X$ and an associated closure operator for equivalent definitions of
merotopic and nearness spaces. In such a way the nearness spaces were originally introduced by Herrlich in \cite{herrlich1} 
and \cite{herrlich2}. For details on the various equivalent formulations we refer the reader to  \cite{herrlich2}.\par
It is also due to Herrlich that $\mathbf{Near}$ is bireflective in $\mathbf{Mer}$. This theorem is usually proved by constructing the 
bireflective modification of a given merotopic space with respect to $\mathbf{Near}$ via transfinite induction (cf. \cite{herrlich2}*{Theorem 8.1} 
or \cite{preuss1}*{Theorem 3.2.2.5}).\par
We want to give a different proof here, which is based on Zorn's Lemma. We begin with an easy lemma.
\begin{lemma*}
If $\mu_1$ and $\mu_2$ are two nearness-structures on the set $X$ then
\begin{equation*}
\mu=\set*{\A\in c(X):\exists \A_1\in \mu_1, \exists \A_2\in \mu_2 \ \mathrm{such\ that} \ \A_1\wedge\A_2\prec\A}
\end{equation*}
is again a nearness-structure on $X$ that contains $\mu_1$ and $\mu_2$. Moreover, every merotopic structure on $X$ 
that contains $\mu_1$ and $\mu_2$ must also contain $\mu$.
\end{lemma*}
Note that if we already knew that $\mathbf{Near}$ is bireflective in $\mathbf{Mer}$ this lemma would be an immediate consequence of the 
general way of constructing initial objects in $\mathbf{Mer}$ (cf. \cite{preuss1}*{Remark 3.2.2.2 \circled{2} and Theorem 3.2.2.1}), but since we want
to use it to show the bireflectivity result we have to give a direct proof, which can be easily done as follows.

\begin{Proof}
Obviously we have $\mu_1, \mu_2\ssq \mu$ and $(X,\mu)$ is easily seen to be a merotopic space. Now pick any $\A\in \mu$. 
By definition there are $\A_1\in \mu_1$ and $\A_2\in \mu_2$ such that $\A_1\wedge\A_2\prec\A$. It follows that $\A_i^{\prime}=
\set*{\interior_{\mu_i}(A):A\in \A_i}\in \mu_i$ for $i=1,2$.\par
Take $A_i\in \A_i$ arbitrarily for $i=1,2$. If $x\in \interior_{\mu_1}(A_1)\cap\interior_{\mu_2}(A_2)$ then $\set*{A_i,X\sm\set{x}}\in \mu_i$
for $i=1,2$ and since $\set*{A_1,X\sm\set{x}}\wedge\set*{A_2,X\sm\set{x}}\prec\set*{A_1\cap A_2,X\sm\set{x}}$ it follows that 
$x\in \interior_{\mu}(A_1\cap A_2)$.\par
Thus we have $\interior_{\mu_1}(A_1)\cap\interior_{\mu_2}(A_2)\ssq \interior_{\mu}(A_1\cap A_2)$ and because of $\A_1\wedge\A_2\prec\A$
we find $A\in \A$ such that $A_1\cap A_2\ssq A$ and hence $\interior_{\mu_1}(A_1)\cap\interior_{\mu_2}(A_2)\ssq \interior_{\mu}(A)$.
So we have that $\A_1^{\prime}\wedge\A_2^{\prime}$ is a refinement of $\set*{\interior_{\mu}(A):A\in \A}$ and hence the latter set 
belongs to $\mu$ which shows that $(X,\mu)$ is indeed a nearness space. The ``moreover'' part is clear.
\end{Proof}

Now we are ready to prove the bireflectivity of $\mathbf{Near}$ in $\mathbf{Mer}$.
\begin{theorem*}
$\mathbf{Near}$ is bireflective in $\mathbf{Mer}$.
\end{theorem*}

\begin{Proof}
Let $(X,\mu)$ be any merotopic space and put 
\begin{equation*}
\mathcal{M}=\set*{\nu \ssq c(X):(X,\nu)\in \mathbf{Near} \ \mathrm{and} \ \nu\ssq \mu}.
\end{equation*}
The set $\mathcal{M}$ is partially ordered by inclusion and $\mathcal{M}$ is non-empty, because $\set*{\A\ssq \mathcal{P}(X):X\in \A}$
is an element of $\mathcal{M}$.\par
If $\mathcal{S}$ is any non-empty chain in $\mathcal{M}$ then it is easy to see that $\bigcup\mathcal{S}$ is again in $\mathcal{M}$ and
hence by Zorn's Lemma there is a maximal element $\tilde{\mu}$ of $\mathcal{M}$.\par
Since $\tilde{\mu}\ssq \mu$ the identity map $\id_X:(X,\mu) \rightarrow (X,\tilde{\mu})$ is uniformly continuous.\par
Now if $(Y,\nu)$ is another merotopic space and $f:(X,\mu) \rightarrow (Y,\nu)$ is uniformly continuous we can put 
\begin{equation*}
\mu_f=\set*{\B\in c(X):\exists \A\in \nu \ \mathrm{with} \ f^{-1}[\A]\prec \B}
\end{equation*}
and show exactly as in the proof from \cite{preuss1}*{Theorem 3.2.2.5} that $(X,\mu_f)$ is a nearness space. Next we define
\begin{equation*}
\bar{\mu}=\set*{\A\in c(X):\exists \A_1\in \tilde{\mu}, \exists \A_2\in \mu_f \ \mathrm{such\ that} \ \A_1\wedge\A_2\prec\A}.
\end{equation*}
By the preceding lemma $(X,\bar{\mu})$ is a nearness space and $\tilde{\mu}, \mu_f\ssq \bar{\mu}$. Since $f$ is uniformly continuous 
with respect to $\mu$ and $\nu$ it follows that $\mu_f\ssq \mu$ and because $\tilde{\mu}$ is also contained in $\mu$ it follows that 
$\bar{\mu}\ssq \mu$, in other words $\bar{\mu}\in \mathcal{M}$ and by the maximality of $\tilde{\mu}$ we must have $\bar{\mu}=\tilde{\mu}$.
Hence $\mu_f\ssq \tilde{\mu}$ which implies that $f:(X,\tilde{\mu}) \rightarrow (Y,\nu)$ is uniformly continuous.
Thus we have shown that $(X,\tilde{\mu})$ is our desired bireflective modification of $(X,\mu)$ with respect to $\mathbf{Near}$.
\end{Proof}

\address
\email

\end{document}